\documentclass{amsart}
\usepackage{latexsym, amssymb}
\usepackage{graphics}

\theoremstyle{plain}
\newtheorem{thm}{Theorem}

\newtheorem{cor}{Corollary}

\newtheorem{quest}{Question}
\theoremstyle{defifinition}
\newtheorem{defn}{Definition}

\begin{document}
\begin{center}
{\bf On the Existence of Finite Type Link Homotopy Invariants}\\
\vspace{.2in}
{\footnotesize BLAKE MELLOR}\\
{\footnotesize Honors College}\\
{\footnotesize Florida Atlantic University}\\
{\footnotesize 5353 Parkside Drive}\\
{\footnotesize Jupiter, FL  33458}\\
{\footnotesize\it  bmellor@fau.edu}\\
\vspace{.1in}
{\footnotesize DYLAN THURSTON}\\
{\footnotesize Department of Mathematics}\\
{\footnotesize Harvard University}\\
{\footnotesize Cambridge, MA  02138}\\
{\footnotesize\it dpt@math.harvard.edu}\\
 
\vspace{1in}
{\footnotesize ABSTRACT}\\
{\ }\\
\parbox{4.5in}{\footnotesize \ \ \ \ \ We show that for links with at most 5 components, the only
finite type homotopy invariants are products of the linking numbers.  In contrast, we show that for
links with at least 9 components, there must exist finite type homotopy invariants which are {\it not}
products of the linking numbers.  This corrects the errors of the first author in \cite{me1, me2}.
\noindent {\it Keywords:}  Finite type invariants; link homotopy.}\\
\end{center}
\tableofcontents
\section{Introduction} \label{S:intro}
In \cite{me1, me2} the first author claimed, erroneously, that there are no finite type link
homotopy or concordance invariants other than the pairwise linking numbers (and their products).  However,
the proofs of this result in both of these paper contained a serious algebraic error.  The purpose of this
paper is to show the opposite - in fact, there {\it do} exist finite type link homotopy (and, hence,
concordance) invariants other than the linking numbers.  However, the proof is not constructive; it
is still an open problem to actually construct such an invariant (see Section~\ref{S:questions}).

There have been many excellent introductions to the theory of finite type invariants, such as
\cite{bl, bn1, cd}; we will not try to replicate them here.  We will provide a few basic definitions in
order to clarify our notation and terminology.  It should be mentioned that our approach and proofs are
combinatorial in nature.

\subsection{Singular Links} \label{SS:singular}
Recall that, in the most general sense, a {\it link invariant} is a map from
the set of equivalence classes of links under isotopy to another set
$G$.  We will need to have some additional structure on $G$.  For our purposes,
$G$ will be the field of complex numbers $\mathbb{C}$.  In this theory, it is also convenient
to look at invariants of {\it regular} isotopy (i.e. links with framing), rather than just
isotopy.  So we will not allow the first Reidemeister move.
We first note that we can extend any link invariant to an invariant of
{\it singular} links, where a singular link is an immersion of several copies of $S^1$ into
3-space which is an embedding except for a finite number of isolated double
points.  Given a link invariant $v$, we extend it via the relation:
$$\includegraphics{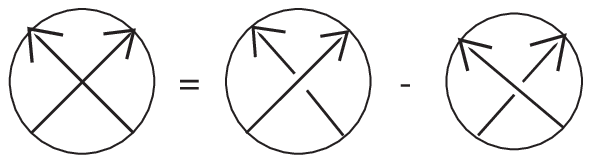}$$
An invariant $v$ of singular links is then said to be of {\it finite type}, if
there is an integer $d$ such that $v$ is zero on any link with more than $d$
double points.  $v$ is then said to be of {\it type} $d$.  We denote by $V_d$
the vector space over $\mathbb{C}$ generated by $\mathbb{C}$-valued finite type invariants of type $d$.  We
can completely understand the space of $\mathbb{C}$-valued finite type invariants by understanding
all of the quotient spaces $V_d/V_{d-1}$.

\subsection{Link homotopy and link concordance}
The idea of {\it link homotopy} (or just {\it homotopy}) was introduced by Milnor~\cite{mi}.
Two links are homotopic if one can be transformed into the other through a sequence of ambient
isotopies of $S^3$ and crossing changes of a component with itself (but {\it not} crossing
changes of different components).  Habegger and Lin~\cite{hl} succeeded in classifying links up to
homotopy.  We construct a theory of finite type invariants in exactly the same way as before; the
difference is that the invariants are trivial when evaluated on a link with a singularity in which a
component intersects itself.  In this case, the two "resolutions" of the singular point are
homotopically equivalent, so the value of a homotopy invariant on their difference is zero.  We
will denote the vector space of type $d$ link homotopy invariants by $V_d^h$.

\begin{defn} \label{D:concordance}
Consider two k-component links $L_0$ and $L_1$.  These can be thought of as embeddings:
$$L_i: \bigsqcup_{i=1}^k S^1 \hookrightarrow {\bf R}^3$$
A {\bf (link) concordance} between $L_0$ and $L_1$ is an embedding:
$$H: \left({\bigsqcup_{i=1}^k S^1}\right) \times I \hookrightarrow {\bf R}^3 \times I$$
such that $H(x,0) = (L_0(x),0)$ and $H(x,1) = (L_1(x),1)$.
A concordance is an isotopy if and only if H is level preserving; i.e. if the image of $H_t$ is
a link at level $t$ for each $t \in I$.
\end{defn}
We will denote the vector space of type $d$ link concordance invariants by $V_d^c$.
 
\subsection{Unitrivalent diagrams} \label{SS:unitrivalent}

It is a marvelous fact that the vector spaces $V_d^*/V_{d-1}^*$ can be given relatively simple
combinatorial descriptions in terms of {\it unitrivalent diagrams}.  These are spaces of
unitrivalent graphs (with colored endpoints and oriented vertices) with various relations imposed upon them.
That these descriptions are
isomorphic to the original vector spaces is largely due to Kontsevich and his integral (see \cite{cd}
for an excellent exposition of the Kontsevich integral).  The description of the space for link homotopy
was developed by Bar-Natan and others \cite{bn1, bgrt}.  The modification for concordance was found by
Habegger and Masbaum \cite{hm}.  For a more detailed development, see \cite{me1, me2}.

\begin{figure}
    $$\includegraphics{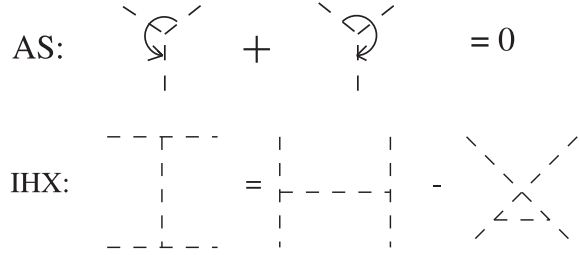}$$
    \caption{AS and IHX relations} \label{F:ihx}
    \end{figure}

\begin{figure}
    $$\includegraphics{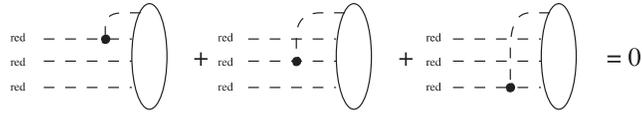}$$
    \caption{The link relation for unitrivalent diagrams} \label{F:link_rel}
    \end{figure}

\begin{defn} \label{D:homotopy}
$B^h$ is defined as the vector space of (disjoint unions of) unitrivalent diagrams modulo the
following relations:
\begin{itemize}
    \item  The antisymmetry (AS) relation (see Figure~\ref{F:ihx}).
    \item  The IHX relation (see Figure~\ref{F:ihx}).
    \item  The link relation (see Figure~\ref{F:link_rel}), where the sum is over {\it all} univalent
    vertices of the diagram with the same color.  Another example of diagrams appearing in a link
    relation can be found in Figure~\ref{F:arise}.
    \item  Any diagram with a loop is trivial.
    \item  Any diagram with a connected component which has two univalent vertices of the same color
    is trivial.
\end{itemize}
The {\bf degree} $d$ of a diagram in $B^h$ is defined to be one half of the number of vertices of the diagram.
Let $B_d^h$ be the vector space of unitrivalent diagrams of degree $d$ (notice that all of the relations
involve diagrams of the same degree, so they apply equally well to $B_d^h$).  So $B^h$ is just the
graded vector space $\bigoplus_{d=1}^{\infty} B_d^h$.  We define $B^c$ to be the space of unitrivalent diagrams
modulo only the first four relations (so components can have multiple endpoints with the same color), and
$B^c = \bigoplus_{d=1}^{\infty} B_d^c$ in the same way.
\end{defn}
\begin{thm} \label{T:homotopy} \cite{bn1, bgrt, hl}
$B_d^h \cong V_d^h/V_{d-1}^h$, and $B_d^c \cong V_d^c/V_{d-1}^c$.
\end{thm}
 
\section{Non-existence results for $B^h$ and $B^c$} \label{S:nonexistence}
Now that we have properly defined the spaces $B^h$ and $B^c$ of unitrivalent diagrams for link homotopy, we
want
to analyze them more closely.  Let $B^h(k)$ (respectively $B^c(k)$) denote the space of unitrivalent
diagrams for link homotopy (resp. concordance)
with $k$ possible colors for the univalent vertices (i.e. we are looking at links with $k$ components).
\subsection{Previous results for $B^h(k)$} \label{SS:3and4}
\begin{figure}
    $$\includegraphics{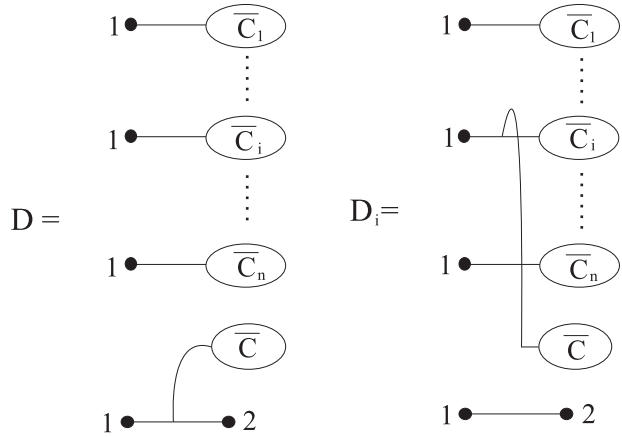}$$
    \caption{Diagrams arising from the link relation} \label{F:arise}
    \end{figure}

Consider a diagram $D \in B^h(k)$.  Each component of $D$
is a tree diagram with at most one endpoint of each color.  Since a unitrivalent tree with $n$
endpoints has $2n-2$ vertices, and hence degree $n-1$, $D$ cannot have any components of degree greater
than $k-1$.
{\bf Notation:}  Before we continue, we will introduce a bit of notation which will be useful in this
section.  Given a unitrivalent diagram $D$, we define $m(D;i,j)$ to be the number of components
of $D$ which are simply line segments with ends colored $i$ and $j$, as shown below:
$$i-----j$$
We call these components {\it struts}.
Recall the following (correct) results from \cite{me1}.  We include the proofs for completeness, and as a
warm up for the more complicated proof in Section~\ref{SS:k=5}:
\begin{thm} \label{T:3comp}
If D has a component C of degree k-1 (with $k \geq 3$), then D is trivial in $B^h(k)$.
\end{thm}
{\sc Proof:}  $C$ has one endpoint of each color $1,2,...,k$.  Without loss of generality, we
may assume that $C$ has a branch as shown, where $\bar{C}$ denotes the remainder of $C$:
$$C:\ \ \begin{matrix} \bar{C} \\ | \\ | \\ 1-----2 \end{matrix}$$
We are going to apply the link relation with the color 1.  Let $\{C_1,...,C_n\}$ be the components of
$D$ with an
endpoint colored 1.  So, ignoring the other components of $D$, we have the diagrams of
Figure~\ref{F:arise} (where $\bar{C_i}$ denotes all of $C_i$ except for the endpoint colored 1).
The link relation then implies that $D + \sum{D_i} = 0$.  If $C_i$ is just a line segment with endpoints
colored 1 and 2, then $D_i = D$.  Otherwise, $\bar{C_i}$ will have an
endpoint of some color $j \in {3,...,k}$.  In this case, since $\bar{C}$ has an endpoint of each color
3,...,$k$, including
$j$, $D_i$ will have a component with two endpoints colored $j$, and hence be trivial in $B^h$.
Therefore, we find that $D+m(D;1,2)D = 0$ where $m(D;1,2) \geq 0$.
We can divide both sides by $1+m(D;1,2)$ (since we are working over a field of characteristic 0) to
conclude that $D = 0$.  $\Box$
 
\begin{thm} \label{T:4comp}
If D has a component C of degree k-2 (with $k \geq 4$), then D is trivial in $B^h(k)$.
\end{thm}
{\sc Proof:}  Without loss of generality, $C$ has endpoints colored $1,2,...,k-1$.  We will
prove the lemma by inducting on $m(D;1,k)$;
inducting among the set of diagrams having a component with endpoints colored $1,2,...,k-1$.
As in the previous theorem, we may assume that $C$ has a branch as shown:
$$C:\ \ \begin{matrix} \bar{C} \\ | \\ | \\ 1-----2 \end{matrix}$$
And conclude that $D+\sum{D_i}=0$, where the $D_i$ are defined as before.  Since $\bar{C}$ contains
endpoints of all colors except 1, 2, and $k$, $D_i$ has two endpoints of the same color (and hence is trivial)
unless $C_i$ has one of the following 3 forms (as in Theorem~\ref{T:3comp}):
$$(1)\ \ C_i = \ 1-----2$$
$$(2)\ \ C_i = \ 1-----k$$
$$(3)\ \ C_i = \ \begin{matrix} k \\ | \\ | \\ 1-----2 \end{matrix}$$
In the first case, $D_i = D$; and in the second case, $D_i = D'$, where $D'$ is the same as $D$ except
that:
\begin{itemize}
    \item $C$ is replaced by a component $C'$ identical to it except that the endpoint colored
2 in $C$ is colored $k$ in $C'$ (so $\bar{C'} = \bar{C}$).
    \item A line segment with endpoints colored 1 and $k$ has been replaced by a line segment
with endpoints colored 1 and 2.  In other words, $m(D';1,2) = m(D;1,2)+1$ and $m(D';1,k) = m(D;1,k)-1$.
\end{itemize}
In the third case, $D_i$ has a component of
degree $k-1$, and so is trivial by the previous theorem.  Therefore, we find
that $D + m(D;1,2)D + m(D;1,k)D' = 0$.  If $m(D;1,k) = 0$ we conclude, as before, that $D$ is trivial
modulo the link relation, which proves the base case of our induction.

\begin{figure}
    $$\includegraphics{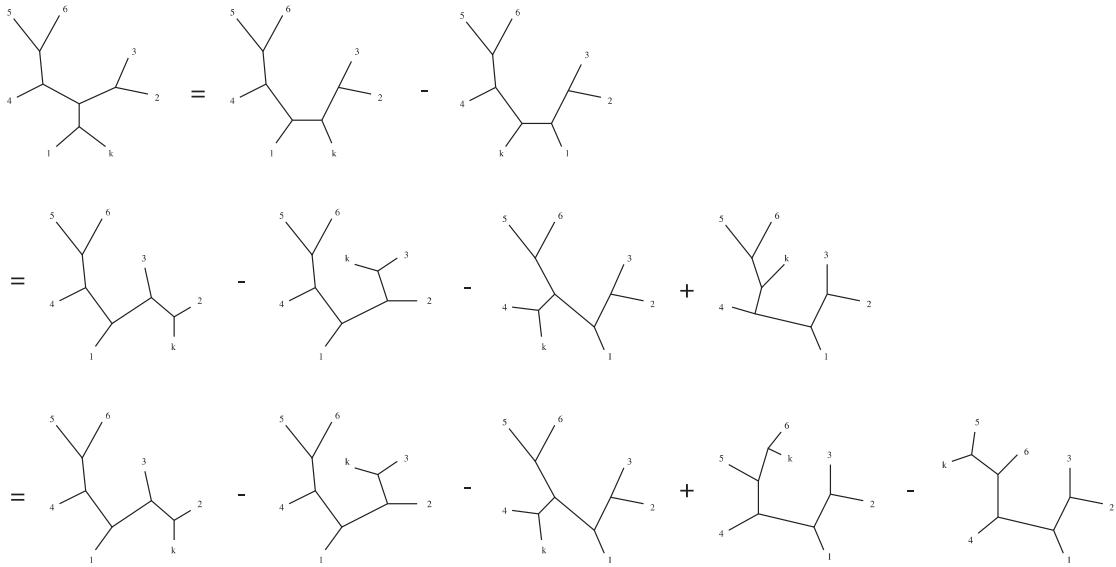}$$
    \caption{Using the IHX relation to decompose a diagram} \label{F:expand}
    \end{figure}

For the inductive step, we use the IHX relation on $C'$ to decompose $D' = \sum_{i\neq 1,2,k}{\pm D_i'}$,
where $D_i'$ is the same as $D'$ except that $C'$ has been replaced by a component $C_i'$ with endpoints
of the same colors (although arranged differently), and a branch as shown:
$$C_i':\ \ \begin{matrix} \bar{C_i'} \\ | \\ | \\ i-----k \end{matrix}$$
(The decomposition is simply a matter of letting the endpoint colored $k$ ``travel'' the tree -
see Figure~\ref{F:expand} for an example.)  In particular, $m(D_i';a,b) = m(D';a,b)$ for all colors
$a$ and $b$.
We now apply the link relation to $D_i'$ using color $i$ (and component $C_i'$).
In this case, the only other components which matter (modulo trivial diagrams) are
ones which look like one of the following:
$$(1)\ \ i-----k$$
$$(2)\ \ i-----2$$
$$(3)\ \ \begin{matrix} 2 \\ | \\ | \\ i-----k \end{matrix}$$
As before, the first case gives $D_i'$ again, the third case is trivial by Theorem~\ref{T:3comp},
and the second case gives a diagram $D_i''$ such that:
\begin{itemize}
    \item $C_i'$ is replaced by a component $C_i''$ identical to it except that the endpoint colored
$k$ in $C_i'$ is colored 2 in $C_i''$ (so $\bar{C_i''} = \bar{C_i'}$).
    \item A line segment with endpoints colored $i$ and 2 has been replaced by a line segment
with endpoints colored $i$ and $k$.  In other words, $m(D_i'';i,k) = m(D_i';i,k)+1$ and
$m(D_i'';i,2) = m(D_i';i,2)-1$.
\end{itemize}
Otherwise, $D_i''$ is the same as $D_i'$; in particular, $m(D_i'';1,k) = m(D_i';1,k) = m(D';1,k) =
m(D;1,k)-1$.
Then the link relation tells us that $D_i' + m(D';i,k)D_i' + m(D';2,i)D_i'' = 0$.  Since $D_i''$
has a component of degree $k-2$ with endpoints colored $1,...,k-1$ (namely, $C_i''$), the inductive
hypothesis implies
that $D_i''$ is trivial.  Therefore, $(1+m(D';i,k))D_i' = 0$, so $D_i'$ is
trivial in $B^h(k)$.  This is true for every $i$, so it immediately follows that $D'$, and hence $D$,
are also trivial in $B^h(k)$.  $\Box$
 
So the largest possible degree of a component of a diagram in $B^h(k)$ is $k-3$ (if $k \geq 4$).  In
particular, this means that if $k$ is 3 or 4, then the largest possible degree of a component of a diagram
in $B^h(k)$ is 1.  It is well-known that the pairwise linking numbers are the only type 1 link homotopy
invariants, and are dual to struts via the isomorphism of Theorem~\ref{T:homotopy}.  Their products are
dual to disjoint unions of struts.
So we have as a corollary:
\begin{cor} \label{C:3and4}
On links with at most 4 components, the only finite type homotopy invariants are the pairwise linking
numbers and their products.
\end{cor}

The obvious question is whether this result will generalize to links with more components.  In the next
section we will show, by a rather involved combinatorial argument, that it {\it does} extend to links with
five components.  However, in Section~\ref{S:existence} we will show that it fails for links with more
than 8 components.

{\sc Remark:}  The error in \cite{me1} (replicated in \cite{me2}) was in the attempt to generalize the
result to all $k$.  On page 785 of \cite{me1}, line 15, the possibility that $c=a$ was neglected.  This
adds another term to the sum, which ends up cancelling everything out.  This was pointed out by Alexander
Merkov \cite{mer}.

\subsection{The case of $B^h(5)$} \label{SS:k=5}
In this section we consider the $B^h(5)$.  We know that no diagram in this space has a component of
degree 3 or more.  So the question is whether a diagram can have a component of degree 2.  Any such
component will be a "Y-component" - i.e. a graph with three (colored) univalent vertices connected
to a single trivalent vertex.

\begin{thm} \label{T:k=5}
If $D \in B^h(5)$ has a component C of degree 2, then D is trivial.
\end{thm}
{\sc Proof:}  This proof is significantly more delicate than that for Theorem~\ref{T:4comp}, involving an
extra level of induction.  Without loss of generality, $C$ has endpoints colored $1,2,3$.
$$C:\ \ \begin{matrix} 3 \\ | \\ | \\ 1-----2 \end{matrix}$$
Our first induction is on $m(D;1,4)+m(D;1,5)$;
inducting among the set of diagrams having a component with endpoints colored $1,2,3$.  We will begin by
proving the base case of this induction.

Let $\{C_1,...,C_n\}$ be the other components of $D$ with an
endpoint colored 1.  Then we apply the link relation as in Figure~\ref{F:arise}.  We will apply the
relation along the color 1, fixing the color 3.  This means that $\bar{C}$ (in Figure~\ref{F:arise}) is
just a single univalent vertex, colored 3.  This will be successively attached to the components
$C_i$ to form the diagrams $D_i$ (in the figure $\bar{C_i}$ denotes all of $C_i$ except for the endpoint
colored 1).  We will refer to this operation as "expanding along 1, fixing 3."

So then $D+\sum{D_i}=0$.  $D_i$ has two endpoints of the same color
unless $C_i$ has one of the following 4 forms:
$$(1)\ \ C_i = \ 1-----2$$
$$(2)\ \ C_i = \ 1-----4$$
$$(3)\ \ C_i = \ 1-----5$$
$$(4)\ \ C_i = \ \begin{matrix} a \\ | \\ | \\ 1-----b \end{matrix}\ \ a,b \in \{2,4,5\}$$
In the first case, $D_i = D$.  In the fourth case, $D_i$ has a component of
degree 3, and so is trivial by Theorem~\ref{T:4comp}.

In the second case, $D_i = D_4$, where $D_4$ is the same as $D$ except that:
\begin{itemize}
    \item $C$ is replaced by a component $C'$ identical to it except that the endpoint colored
2 in $C$ is colored 4 in $C'$.
$$C':\ \ \begin{matrix} 3 \\ | \\ | \\ 1-----4 \end{matrix}$$
    \item A strut with endpoints colored 1 and 4 has been replaced by a strut
with endpoints colored 1 and 2.  In other words, $m(D';1,2) = m(D;1,2)+1$ and $m(D';1,4) = m(D;1,4)-1$.
\end{itemize}

{\sc Notation:}  For the remainder of this proof, we will represent diagrams by giving the changes made
from $D$.  We will draw the new component of degree 2 which has replaced $C$ (we will see that for all of
our diagrams, any other components of degree 2 remain unchanged).  Although the total number of
struts is always preserved, some struts have been replaced by others.  We represent a strut by the (unordered)
pair of the colors of its endpoints, and use an arrow to show how the struts have been traded.  For example,
we will represent $D_4$ as follows:
$$D_4 = \begin{matrix} 3 \\ | \\ | \\ 1-----4 \end{matrix} (1,4) \rightarrow (1,2)$$

Finally, in the third case, $D_i = D_5$, which is defined similarly to $D_4$.
$$D_5 = \begin{matrix} 3 \\ | \\ | \\ 1-----5 \end{matrix} (1,5) \rightarrow (1,2)$$
Therefore, we find
that $D + m(D;1,2)D + m(D;1,4)D_4 + m(D;1,5)D_5 = 0$.  If $m(D;1,4)+m(D;1,5) = 0$, then we conclude
$D + m(D;1,2)D = 0$ (since $m(D;i,j) \geq 0$).  Since $m(D;1,2) \geq 0$, we can divide by $1 + m(D;1,2)$
to conclude $D=0$, which proves the base case of our first induction.

Henceforth, for convenience, we will let $m(i,j) = m(D;i,j)$.

We will now assume the inductive hypothesis that any diagram $E \in B^h(5)$
with a component of degree 2 with endpoints colored 1, 2 and 3, and such that
$m(E;1,4)+m(E;1,5) < m(1,4)+m(1,5)$, is trivial.  Our inductive step consists of using this hypothesis
to prove that $D_4$ and $D_5$ are trivial.  This will immediately imply that $D+m(1,2)D=0$, and hence that
$D=0$.

We will prove that $D_4$ is trivial.  The proof that $D_5$ is trivial is very similar.  This proof involves
a second induction.  We will be looking at diagrams which do not have a component with endpoints colored
1, 2 and 3, so are not directly trivial by the (first) inductive hypothesis.  However, we will find that
(modulo the inductive hypothesis), we can effectively "swap" struts in these diagrams so that the
number of struts with certain colors on their endpoints always decreases.  Since there are only a finite
number of such struts, the supply eventually disappears, and we are able to conclude that the diagrams
are trivial.

We begin with $D_4$.  We expand along 3, fixing 1.
$$D_4 + m(3,4)D_4 + m(2,3)D_{42} + m(3,5)D_{45} = 0$$
$$D_{42} = \begin{matrix} 3 \\ | \\ | \\ 1-----2 \end{matrix}\
        \begin{matrix} (1,4) \rightarrow (1,2) \\ (2,3) \rightarrow (3,4) \end{matrix}$$
$$D_{45} = \begin{matrix} 3 \\ | \\ | \\ 1-----5 \end{matrix}\
        \begin{matrix} (1,4) \rightarrow (1,2) \\ (3,5) \rightarrow (3,4) \end{matrix}$$
$D_{42}$ has a component of degree 2 with endpoints colored 1, 2, and 3.  Also,
$m(D_{42};1,4) + m(D_{42};1,5) = (m(1,4)-1) + m(1,5) = m(1,4)+m(1,5)-1$.  So by the inductive
hypothesis, $D_{42} = 0$.  Therefore:
$$D_4 = \frac{-m(3,5)}{1+m(3,4)} D_{45}$$

Consider $D_{45}$.  We expand along 3, fixing 5.
$$D_{45} + m(1,3)D_{45} + m(2,3)D_{452} + m(3,4)D_{454} = 0$$
$$D_{452} = \begin{matrix} 3 \\ | \\ | \\ 2-----5 \end{matrix}\
        \begin{matrix} (1,4) \rightarrow (1,2) \\ (3,5) \rightarrow (3,4) \\
        (2,3) \rightarrow (1,3) \end{matrix}$$
$$D_{454} = \begin{matrix} 3 \\ | \\ | \\ 4-----5 \end{matrix}\
        \begin{matrix} (1,4) \rightarrow (1,2) \\ (3,5) \rightarrow (3,4) \\
        (3,4) \rightarrow (1,3) \end{matrix} \Rightarrow \begin{matrix} (1,4) \rightarrow (1,2) \\
        (3,5) \rightarrow (1,3) \end{matrix}$$
Therefore:
$$D_4 = \frac{-m(3,5)}{1+m(3,4)}\frac{1}{1+m(1,3)}\left({-m(2,3)D_{452}-m(3,4)D_{454}}\right)$$

Neither of the new diagrams are trivial inductively, so we will need to analyze both of them.  First,
we consider $D_{452}$.  We will show that, modulo the inductive hypothesis, we can swap a strut (3,5) (i.e.
a strut with endpoints colored 3 and 5) for a strut (3,4), while simultaneously swapping a strut (2,4) for a
strut (2,5).  We begin by expanding along 2, fixing 3.
$$D_{452} + m(2,5)D_{452} + (m(1,2)+1)D_{4521} + m(2,4)D_{4524} = 0$$
$$D_{4521} = \begin{matrix} 3 \\ | \\ | \\ 2-----1 \end{matrix}\
        \begin{matrix} (1,4) \rightarrow (1,2) \\ (3,5) \rightarrow (3,4) \\ (2,3) \rightarrow (1,3) \\
        (1,2) \rightarrow (2,5) \end{matrix} \Rightarrow \begin{matrix} (1,4) \rightarrow (2,5) \\
        (3,5) \rightarrow (3,4) \\ (2,3) \rightarrow (1,3) \end{matrix}$$
$$D_{4524} = \begin{matrix} 3 \\ | \\ | \\ 2-----4 \end{matrix}\
        \begin{matrix} (1,4) \rightarrow (1,2) \\ (3,5) \rightarrow (3,4) \\
        (2,3) \rightarrow (1,3) \\ (2,4) \rightarrow (2,5) \end{matrix}$$
By the inductive hypothesis, $D_{4521} = 0$.  Therefore:
$$D_{452} = \frac{-m(2,4)}{1+m(2,5)}D_{4524}$$

Now consider $D_{4524}$.  We expand along 3, fixing 2.
$$D_{4524} + (m(3,4)+1)D_{4524} + (m(1,3)+1)D_{45241} + (m(3,5)-1)D_{45245} = 0$$
$$D_{45241} = \begin{matrix} 3 \\ | \\ | \\ 2-----1 \end{matrix}\
        \begin{matrix} (1,4) \rightarrow (1,2) \\ (3,5) \rightarrow (3,4) \\ (2,3) \rightarrow (1,3) \\
        (2,4) \rightarrow (2,5) \\ (3,4) \rightarrow (1,3) \end{matrix} \Rightarrow
        \begin{matrix} (1,4) \rightarrow (1,2) \\ (3,5) \rightarrow (1,3) \\ (2,3) \rightarrow (1,3) \\
        (2,4) \rightarrow (2,5) \end{matrix}$$
$$D_{45245} = \begin{matrix} 3 \\ | \\ | \\ 2-----5 \end{matrix}\
        \begin{matrix} (1,4) \rightarrow (1,2) \\ (3,5) \rightarrow (3,4) \\
        (2,3) \rightarrow (1,3) \\ (2,4) \rightarrow (2,5) \\ (3,5) \rightarrow (3,4) \end{matrix}$$
By the inductive hypothesis, $D_{45241} = 0$.  Therefore:
$$D_{4524} = \frac{-(m(3,5)-1)}{2+m(3,4)}D_{45245}$$
Combining this with the result of the previous step, we have:
$$D_{452} = \frac{-m(2,4)}{1+m(2,5)}\frac{-(m(3,5)-1)}{2+m(3,4)}D_{45245}
    = \frac{m(2,4)}{1+m(2,5)}\frac{(m(3,5)-1)}{2+m(3,4)}D_{45245}$$

Notice that the component of degree 2 in $D_{45245}$ is the same as that in $D_{452}$, so we could start the
whole procedure again.  Also notice that we have swapped the struts as we wanted.  Inductively, we can see
that:
$$D_{452} = \prod_{k=1}^n{\frac{(m(2,4)-k+1)}{k+m(2,5)}\frac{(m(3,5)-k)}{1+k+m(3,4)}}D_{452(45)^n}$$
Eventually, $n > m(3,5)$, so $D_{452} = 0$.  This means that:
$$D_4 = \frac{-m(3,5)}{1+m(3,4)}\frac{-m(3,4)}{1+m(1,3)}D_{454}$$

Now we examine $D_{454}$.  We expand along 4, fixing 3.
$$D_{454} + m(4,5)D_{454} + (m(1,4)-1)D_{4541} + m(2,4)D_{4542} = 0$$
$$D_{4541} = \begin{matrix} 3 \\ | \\ | \\ 4-----1 \end{matrix}\
        \begin{matrix} (1,4) \rightarrow (1,2) \\ (3,5) \rightarrow (1,3) \\
        (1,4) \rightarrow (4,5) \end{matrix}$$
$$D_{4542} = \begin{matrix} 3 \\ | \\ | \\ 4-----2 \end{matrix}\
        \begin{matrix} (1,4) \rightarrow (1,2) \\ (3,5) \rightarrow (1,3) \\
        (2,4) \rightarrow (4,5) \end{matrix}$$
Notice that $D_{4542}$ has a degree 2 component which is the same, up to sign, as $D_{4524}$.  So, by using
the same kind of argument used for $D_{452}$, we can show that $D_{4542} = 0$ (we will look at the
diagrams $D_{4542(54)^n}$).  So:
$$D_{454} = \frac{m(1,4)-1}{1+m(4,5)}(-D_{4541})$$

Therefore, we can express $D_4$ in terms of $D_{4541}$:
$$D_4 = \frac{m(3,5)}{1+m(3,4)}\frac{m(3,4)}{1+m(1,3)}\frac{m(1,4)-1}{1+m(4,5)}(-D_{4541})$$
Notice that we have swapped a strut (1,4) for a strut (4,5), and a strut (3,5) for a strut (3,4).  Since
$D_{4541}$ has the same degree 2 component (up to sign) as $D_4$, we can repeat the same sequence
of operations.  Inductively, we will find that:
$$D_4 = \prod_{k=1}^n{\frac{(m(3,5)-k+1)}{1+m(3,4)}\frac{m(3,4)}{k+m(1,3)}\frac{(m(1,4)-k)}{k+m(4,5)}}
(-1)^nD_{4(541)^n}$$
Eventually, $n > m(1,4)$, so $D_4 = 0$.

By a similar argument, $D_5$ will equal 0, so we can conclude that $D + m(1,2)D=0$.  Therefore, $D=0$,
which completes the proof.  $\Box$

\begin{cor} \label{5comp}
On links with at most 5 components the only finite type invariants are the pairwise linking numbers and
their products.
\end{cor}

It appears to be difficult to extend the proof of Theorem~\ref{T:k=5} to $B^h(6)$.  In this case, the IHX
relation comes into play, and it is not clear that one can decrease $m(1,4)+m(1,5)$ monotonically.
\subsection{Previous results for $B^c(k)$} \label{SS:2and3}
Unlike for homotopy, there is no {\it a priori} limit on the size of the components of diagrams in $B^c(k)$.
However, we are able to prove some non-existence results for small values of $k$.  The proofs of the following
results can be found in \cite{me2}:
\begin{thm} \label{T:123}
The only nontrivial diagrams in $B^c(k)$ for $k=1,2,3$ are disjoint unions of struts.  In other words,
any diagram with a component of degree greater than 1 is trivial.
\end{thm}

As with homotopy, it seems difficult to extend the methods used to prove Theorem~\ref{T:123} to higher values
of $k$.  The attempt in \cite{me2} fell prey to the same error as in \cite{me1} (see the Remark in
Section~\ref{SS:3and4}).

{\sc Remark:}  If we allow the first Reidemeister move, we can prove that $B^c(1)$ (with the new relation)
is trivial, confirming a result of \cite{ng} that the Arf invariant (which is $\mathbb{Z}_2$-valued) is the
only finite type knot concordance invariant.
 
\section{Existence results for $B^h$} \label{S:existence}

In this section we demonstrate the existence of non-trivial diagrams in $B^h(k)$ which are {\it not} just
the products of small components, for $k \geq 9$.  Since $B^h$ is just a quotient of $B^c$, this implies
the existence of non-trivial diagrams in $B^c$.  This can also be proved directly using the same methods, but
we will leave that as an exercise for the reader.  The arguments used are not constructive; we simply
use a counting argument to show that (within a certain subspace) there are more diagrams than relations.
Since all the relations are just linear combinations, we have a homogeneous system of linear equations with
more equations than unknowns, and we conclude that there must be non-trivial solutions.  Each such solution
corresponds to some finite type invariant which is not just a product of linking numbers.

We consider the subspace $Y^h(k)$ of $B^h(k)$ which is spanned by diagrams which have a single Y-component
(degree 2 component) and all other components are struts (degree 1).  This is the space of diagrams which
have exactly one trivalent vertex.  Since all the relations of $B^h$ preserve the number of trivalent
vertices (i.e. any two diagrams in a given relation have the same number of trivalent vertices), $Y^h(k)$
is closed under the relations of $B^h$.  We will show that $Y^h(k)$ contains non-trivial diagrams for
$k \geq 9$.

We count the number of diagrams in $Y^h(k)$ which have exactly $n$ struts.  We count these diagrams by
counting the number of ways of coloring the endpoints of the Y-component (i.e. of choosing 3 distinct colors),
and then counting the number of ways of choosing the $n$ struts (i.e. of choosing $n$ pairs of distinct
colors).

Notice that this count does not distinguish the orientation of the trivalent vertex.
This would double the number of diagrams, except that the new ones are just the negatives of the previous
ones by the antisymmetry relation.  So we will leave them out, and simply not count the antisymmetry
relations among our relations.  Since our diagrams only have one trivalent vertex, there are no
IHX relations.  Also, since our diagrams have no loops, and the endpoints of any component are given distinct
colors, we can ignore the fourth and fifth relations of Definition~\ref{D:homotopy}.  This means that the
only relations we need to count are the link relations.

Let $u(n,k)$ be the number of elements of $Y^h(k)$ which have $n$ struts (i.e. our number of "unknowns").
There are $\binom{k}{3}$ ways of choosing the labels for the Y-component.  There are $\binom{k}{2}$ possible
struts.  The number of ways of selecting $n$ of them, with repetition allowed, is simply
$\binom{\binom{k}{2} + n -1}{n}$.  Putting these together, we find:
$$u(n,k) = \binom{k}{3}\binom{\binom{k}{2} + n -1}{n}$$

Now we want to count the link relations among these diagrams.  Notice that the diagrams in $Y^h(k)$ with
$n$ struts are exactly the diagrams in $Y^h(k)$ with $2n+3$ endpoints.  Since the link relation
preserves the number
of endpoints, as well as the number of univalent vertices, if one diagram in a relation is an element of
$Y^h(k)$ with $n$ struts, so is every other diagram in the relation.  Let $r(n,k)$ be the number of link
relations among elements in $Y^h(k)$ which have $n$ struts.  We can think of one of these relations as
consisting of $n+1$ struts, together with one "special" strut.
The special strut will have a distinguished endpoint.  The link relation is created by attaching this
endpoint in turn to all the other struts which have an endpoint of the same color, forming a series of
diagrams with a single Y-component.  An example is shown below:
$$\begin{matrix}
    3-----2* \\ \\ 2-----1 \\ \\ 2-----4 \\ \\ 5-----6
    \end{matrix} \ \ \longrightarrow \ \ \begin{matrix}
        \begin{matrix}
            3 \\ | \\ | \\ 2-----1
            \end{matrix} \\ \\
        2-----4 \\ \\
        5-----6
        \end{matrix} \ \ + \ \ \begin{matrix}
        \begin{matrix}
            3 \\ | \\ | \\ 2-----4
            \end{matrix} \\ \\
        2-----1 \\ \\
        5-----6
        \end{matrix} \ \ = \ \ 0$$
There are $k(k-1)$ ways of coloring the "special" strut (not
$\binom{k}{2}$, since the endpoints are distinguished).  Then, as before, there are
$\binom{\binom{k}{2} + n}{n+1}$ ways of choosing the other $n+1$ struts.  We conclude that:
$$r(n,k) = k(k-1)\binom{\binom{k}{2} + n}{n+1}$$

To compare these two counts, we look at the quotient of the number of relations by the number of diagrams:
$$\frac{r(n,k)}{u(n,k)} =
\frac{k(k-1)\binom{\binom{k}{2} + n}{n+1}}{\binom{k}{3}\binom{\binom{k}{2} + n -1}{n}} =
\frac{6}{k-2}\frac{\binom{k}{2}+n}{n+1}$$
For a fixed value of $k$, we can look at the limit of this ratio as $n \to \infty$:
$$\lim_{n \to \infty} \frac{r(n,k)}{u(n,k)} = \frac{6}{k-2}$$
If $k \geq 9$, then this limit is less than 1, which means there are more relations than diagrams, so there
must be non-trivial diagrams.  If we plug in $k=9$ and solve for the ratio to be 1, we obtain $n=209$.
So if we have 210 struts (i.e. diagrams of degree 212) there will definitely be nontrivial diagrams.

\begin{thm}
There is a non-trivial homotopy invariant on links with 9 components, of type 212, which is not a product
of linking numbers.
\end{thm}

{\sc Remark:}  In fact, we have slightly overcounted the relations.  We have counted diagrams where the
distinguished endpoint of the special strut has a color which does not appear elsewhere in the diagram, so
it cannot be attached to any other strut to form a Y-component.  However, if we take the limit of the number
of these extra relations divided by $u(n,k)$ as $n$ tends to $\infty$, we get 0.  So removing these relations
from the count does not significantly improve our result.

In general, of course, many of the relations are dependent.  So we would expect that there are also
non-trivial diagrams when $k=8$, and the ratio tends to 1, and possibly for even lower values of $k$.

\section{Questions} \label{S:questions}

\begin{quest}
What is an explicit example of a non-trivial finite type link homotopy invariant which is not a
product of linking numbers?
\end{quest}

Any such invariant would immediately give a finite type invariant
for string links.  Bar-Natan \cite{bn2} has shown that the finite type invariants for string links
are exactly the Milnor invariants, which classify string links up to homotopy.  However, Milnor's invariants,
other than the linking numbers, have indeterminacies which prevent them from being lifted to links as
integer- (or $\mathbb{C}$-) valued invariants.  Apparently, it is possible to find some product in
which the indeterminacies "cancel" and the product can be lifted, which is unexpected.

\begin{quest}
What is the first value of $k$ for which $B^h(k)$ has non-trivial diagrams which are not disjoint unions of
struts?
\end{quest}

We know that such diagrams exist for $k \geq 9$, and that they do {\it not} exist for $k \leq 5$, but the
situation for $k = 6,7,8$ is still unknown.  It seems likely there are non-trivial diagrams in $B^h(8)$,
but as yet we do not have a proof.

\begin{quest}
What is the first value of $k$ for which $B^c(k)$ has non-trivial diagrams which are not disjoint unions of
struts?
\end{quest}

As for $B^h(k)$, we know that such diagrams exist for $k \geq 9$, and do not exist for $k \leq 3$, but the
situation for $4 \leq k \leq 8$ is unknown.

\begin{quest}
Can we refine the methods of Section~\ref{S:existence} to prove that there are nontrivial diagrams for
lower values of $k$?
\end{quest}

We mentioned in Section~\ref{S:existence} that the link relations preserve the number of trivalent vertices
of the diagram.  They also preserve the number of univalent vertices of each color.  Perhaps this could be
used to find smaller "closed" subspaces with fewer dependent relations, so that the ratio of relations to
diagrams is smaller in the limit.

\section{Acknowledgements}
The first author wishes to thank Alexander Merkov for pointing out the error in his earlier papers, and
for his comments on the proof of Theorem~\ref{T:k=5}.

\end{document}